# Addressing Data Engineering Challenges in the Cold-Chain Sector to Reduce its Environmental Impact Through Analytics

**Mathilde MARCY*[(a)], Camille FERTEL[(b)], Thomas SUQUET[(b)], Gérald CAVALIER[(b)]**

[(a)] CEMAFROID; [(b)] TECNEA
Fresnes, 94260, France
*Corresponding author: mathilde.marcy@cemafroid.fr

## Abstract

Data analytics offer a great opportunity for organizations within the cold-chain industry to become more sustainable and reduce their environmental footprint. Yet, despite recent progress in data sciences, they still face challenges to ensure data quality and streamline data processes through adapted full-scale ecosystems, limiting their ability to perform sustainable analytics.

Most data techniques and tools stem from the data-science community and are not always suitably adapted to domain specificities. For instance, nearly all data-quality frameworks are designed to individually consider entities, i.e. representations of real-world elements. However, achieving environmental impact reduction in industrial settings, particularly within the cold-chain sector, often requires modelling and studying entire systems comprised of multiple interdependent entities.

This publication challenges the application of conventional data science and engineering concepts in data ecosystems within the cold-chain sector and proposes an alternative approach to sustainably address prevalent challenges and thus optimize analytical processes.

## 1. Introduction

The cold chain is essential to our modern lives, guaranteeing safe access to medication and perishable food, but unfortunately, it has positive and negative impacts on the environment. For instance, domestic refrigerators electricity consumption account for more than 4% of global world electricity consumption but they also contribute to reducing food losses. In the collective effort to mitigate climate change, the cold chain sector must adopt strategies and technologies to reduce its environmental footprint and even make a positive contribution to environment. Data science and analytics have been identified as a great vector of sustainability for the sector (Bottani et al., 2022; Chaudhuri et al., 2018; Mohan and Amin, 2025), and the reduction potential they offer is quickly increasing with the rapid evolution of data science's technologies and the advent of industry 4.0 (Bhagat et al., 2022; Sartal et al., 2020; Vrchota et al., 2020), already spreading to the cold chain sector (Fatorachian and Pawar, 2025).

Machine learning algorithms can be used to predict technologies' breakdowns (Kale, 2022; Lorenc et al. 2021) and refrigerant leaks (Mtibaa, 2025), to reduce consumption (Zhao et al., 2020), and avoid food losses (Khanna et al., 2025). Internet of things (IoT), coupled with cloud computing and cyber physical systems (CPS), allow real-time monitoring and actions to prevent failure or spoilage (Lorenc et al., 2021; Kavididevi et al., 2024; Gillespie et al., 2023). Digital twins considerably improve experimental capacity thanks to low running cost simulation to optimize processes and resource allocation (Wu et al., 2023).

Democratization of data science throughout the cold chain offers valuable potential to reduce its environmental footprint. However, to be sustainable, data science operations rely on complex data ecosystems (i.e. a network of tools, processes, and rules, that manage the entire data lifecycle, from data sources to final products consumption). Their deployment and use can be challenging and carries a non-negligeable environmental cost (Corbett, 2018; Li and Huang, 2023; Whitehead et al., 2014) . The role of data engineering is to design, operate, and maintain sustainable and performant data ecosystems.

Beyond the inherent complexity of data engineering, the cold chain sector imposes domain-specific requirements that further complexify the design and operations of these data ecosystems. First, organizations must comply with multiple regulatory and standard frameworks, such as the United Nations ATP[1] agreement (ECE/TRANS/347), the European F-Gas regulation (UE/2024/573), and European Pressure Equipment Directive (2014/68/UE), which may impose partially conflicting data management requirements. Second, the potentially severe consequences of technology failure or erroneous decisions, ranging from risks to public health to environmental impacts due to product spoilage, necessitate high data quality and low-latency data processing, particularly for time-sensitive decisions based on real-time monitoring. Third, many cold chain organizations operate across multiple complementary activities, requiring the integration of disparate data sources to accurately assess the environmental impact of their operations. However, the deployment and operation of complex ecosystems capable of accommodating these requirements are often computationally intensive and entail environmental costs related to large-scale data storage (Whitehead et al., 2014), sustained computational workload, and hardware production and disposal, which may ultimately offset an organization's sustainability gains from analytics. To keep the balance in favor of reducing the environmental footprint, data ecosystems must therefore be carefully optimized for the specific requirements of the cold chain.

Despite significant technological advances in data science, many data management techniques and tools remain largely domain-agnostic, complicating their application to the cold chain. For instance, data modelling and data quality methods, which are essential to any efficient data ecosystem, are typically grounded in the notion of individual entities (i.e. real-world elements that can be uniquely identified and represented by single data records) and are designed to operate at this level. Duplicate detection, for example, usually focuses on identifying records that represent the same entity (Christen, 2012; Draisbach et al., 2020). However, many systems of interest in the cold chain are composed of multiple interdependent entities, which may themselves be hierarchical or composite in nature. Capturing this complexity requires fine-grained modelling that many existing techniques do not adequately support, thereby increasing process complexity, computational cost, and, in turn, the environmental footprint of data operations.

The authors argue for the emergence of a *cold chain data management* discipline to enable the sector to fully leverage analytics and modern data technologies for environmental footprint reduction. Such discipline would support the effective adaptation of existing methods and foster the development of solutions aligned with the sector's specific constraints. This paper contributes by identifying domain-specific challenges that hinder the deployment of data science and engineering in the cold chain and by laying initial foundations for the *cold chain data management*.

Section 2 reviews existing research on environmental impact reduction in the cold chain through data-driven approaches. Sections 3 and 4 examine foundational data concepts and their limitations with respect to domain-specific constraints. Section 5 discusses key principles of the proposed discipline, and Section 6 concludes.

## 2. Data Science and Analytics in the Cold-Chain

There is growing interest in the use of analytics and data science to optimize cold chain operations and reduce their environmental impact. Early work has relied on classical modelling approaches; for example, Fabris et al. (2022) and Yao et al. (2023) developed dynamic simulation models to assess the environmental performance, fuel consumption, and emissions of refrigerated vehicles. Predictive models have also been widely applied, particularly to anticipate spoilage and reduce emissions. Kale (2022) used models such as random forests, decision trees, and k-nearest neighbors to predict food quality and detect cold chain disruptions using real-time temperature data. Similarly, Zhao et al. (2020) used ant colony optimization to solve multi-objective routing problems, reducing both distribution costs and carbon emissions. Capo et al. (2020) propose a random forest model to predict the evolution of the K coefficient of refrigerated truck isothermal bodies over time and identifying key ageing factors relevant to emissions reduction.

[1] Agreement on the International Carriage of Perishable Foodstuffs and on the Special Equipment to be Used for such Carriage as amended on 22 June 2024

More advanced technologies are increasingly explored, including blockchain for traceability (Arora et al., 2024) and digital twins for incident detection and localization (Wu et al., 2023). In parallel, machine learning, reinforcement learning, and deep learning approaches are gaining prominence for equipment fault prediction and operational optimization based on real-time IoT data such as temperature readings and hygrometry. Kavididevi et al. (2024) proposed reinforcement learning to predict adverse events and trigger corrective actions, while Gillespie et al. (2023) and Lorenc et al. (2021) leveraged IoT data and neural networks to detect equipment failures and transportation disruptions. Khanna et al. (2025) further combined reinforcement learning, blockchain, and generative AI to significantly reduce spoilage, energy consumption, and emissions.

Most existing studies on the use of analytics to reduce the cold chain's footprint typically examine a single technology, process, or application. As a result, they often overlook the broader data value chain behind. Sustainably leveraging operational data through analytics requires automating and optimizing data processes, which can only be achieved within a data ecosystem, as illustrated in Figure 1.

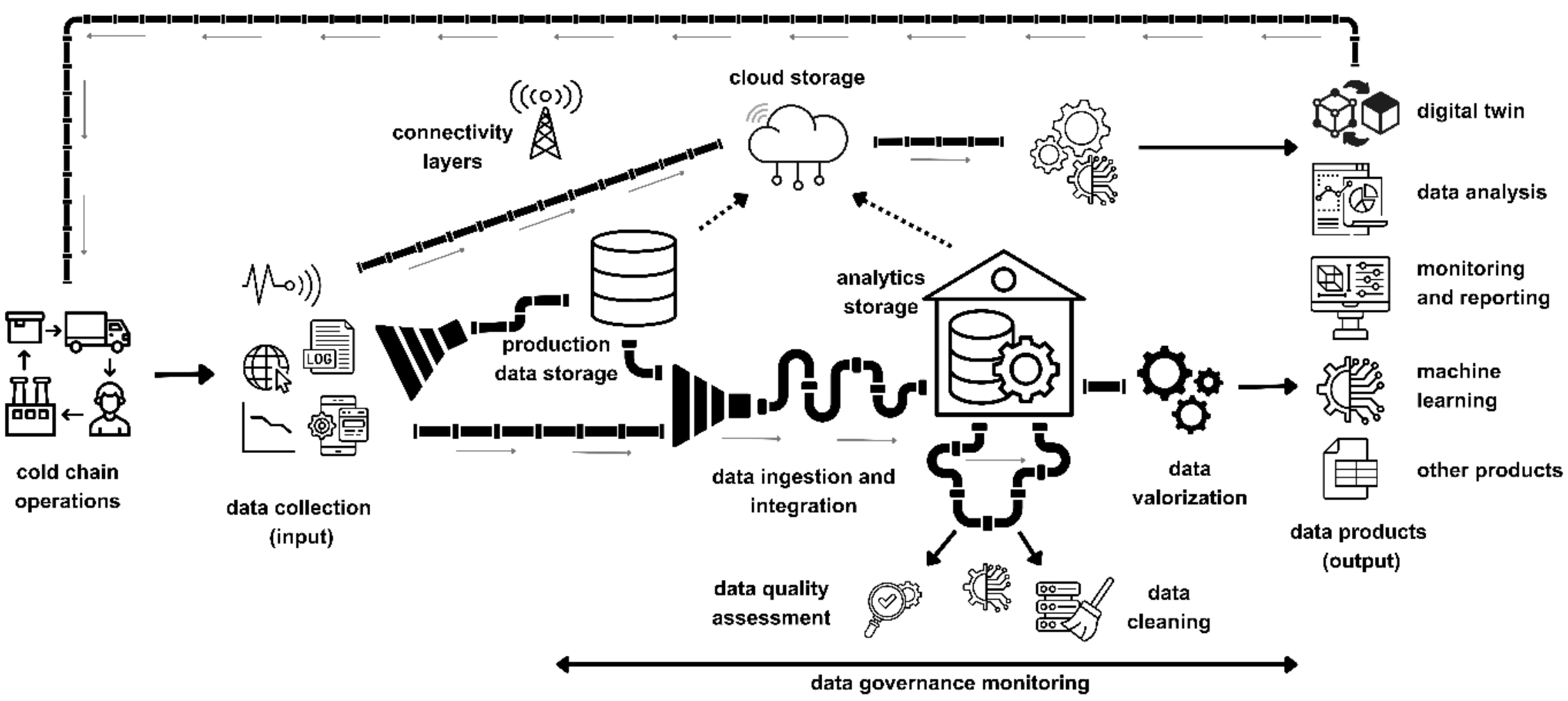

**Figure 1: Typical data ecosystem in the cold chain**

These data ecosystems and their components are essential to the deployment of data products, and themselves carry a significant environmental footprint; they should therefore not be overlooked when assessing the environmental impact of data-driven approaches in the cold chain. The design, development, and maintenance of such ecosystems fall within the scope of data engineering, a discipline grounded in foundational data principles that govern how data is structured, processed, and managed.

## 3. Data Science and Engineering Fundamentals

Although data science and engineering tools evolve at an accelerated pace, a set of foundational concepts continues to underpin most data-driven techniques and processes, including the abstract concept of "entity" which govern most data modelling and quality paradigms. An entity is a representation of a real-world element that can be uniquely identified and is described by attributes, and an entity class is a collection of comparable entities sharing the same attributes and that can among themselves be uniquely identified. Entity classes are inherently connected to data semantics, which provide meaning and context to data. For example, the same data is represented in Tables 1, 2 and 3. However, data from Table 3 is much easier to interpret because its model semantically indicates what real-world elements they represent.

**Table 1. columnar data without semantics**

| A |
|---|
| xy |
| TRUE |
| 10 |
| FC |
| acme |
| zy |
| FALSE |
| 11 |
| FFT |
| Acme |
| xy |
| 11 |
| ACME inc. |
| ab |
| TRUE |
| -21 |
| acme |

**Table 2. data without semantics**

| A | B | C | D | E | F |
|---|---|---|---|---|---|
| C1 | xy | TRUE | 10 | FC | acme |
| C2 | zy | FALSE | 11 | FFT | Acme |
| C3 | xy |  | 11 | F.C. | ACME inc. |
| C4 | ab | TRUE | -21 | fft | acme |

**Table 3. data described with semantics**

| cell_id | serial_number | operational | surface | brand | manufacturer |
|---|---|---|---|---|---|
| C1 | xy | TRUE | 10 | FC | acme |
| C2 | zy | FALSE | 11 | FFT | Acme |
| C3 | xy |  | 11 | F.C. | ACME inc. |
| C4 | ab | TRUE | -21 | fft | acme |

### 3.1. Entity and entity class formalization

In the context of refrigerated transportation, consider an entity class $C$ representing insulated bodies, described by a set of attributes $R_c$.
$R_c = \{serial_number, brand, manufacturing_date, interior_surface, isothermal_class, \ldots\}$
and $c \in C$ an entity representing a specific real-world insulated body, characterized by a data value $c[A]$ for each attribute $A \in R_C$.

For example, from Table 3, if $c[cell_{id}] = C1$ then $c[serial_number] = xy$ and $c[brand] = FC$.

Each entity is uniquely identifiable by a *natural key*, defined as the attribute or set of attributes, that uniquely identifies entities within an entity class. For class $C$, the natural key is $N_C = \{brand, serial_number\}$.

Thus, for any $c_1, c_2 \in C$, if $c_1[N_C] = c_2[N_C]$ (or $c_1[N_C] \approx c_2[N_C]$ to account for data quality issues), then $c_1$ and $c_2$ are duplicates, representing the same real-world insulated body, as two distinct cells cannot share the same brand and serial number.

Most data engineering processes and algorithms are semantic-agnostic, operating on data values regardless of their real-world meaning. From a technical standpoint, training a machine learning model on data from Table 2 or 3 follows the same computational logic, independent of what the data represents. The notion of an entity is not intrinsically embedded in mathematical models or data processing pipelines, which solely do what they are programmed to do.

However, semantics and entity-centric data structure play a critical role upstream in the data ecosystem, particularly in data modelling and data quality management, where data is evaluated in comparison to real-world objects and processes. When semantic structure and data quality are neglected, downstream tools process flawed inputs indiscriminately, leading to potentially inaccurate outcome, a phenomenon commonly known as the "garbage in, garbage out" principle.

### 3.2. Data modelling

Data modelling is the process of defining the structure, semantics, relationships, and constraints of data within a given domain, which serves as a blueprint for all data processes within an ecosystem. Most modelling paradigms rely on the concept of entity; In conceptual modelling, domain experts represent key business entity classes and define them along with their attributes and quality and integrity rules. Figure 2 represents a partial entity-relationship diagram for refrigerated vehicles.

Take the example of Table 1. Although it includes the same data as Tables 2 and 3, it is doubtful that any model can be learned from this data, nor that its quality could be automatically evaluated, since without logical structure, patterns within the data are not easily recognizable.

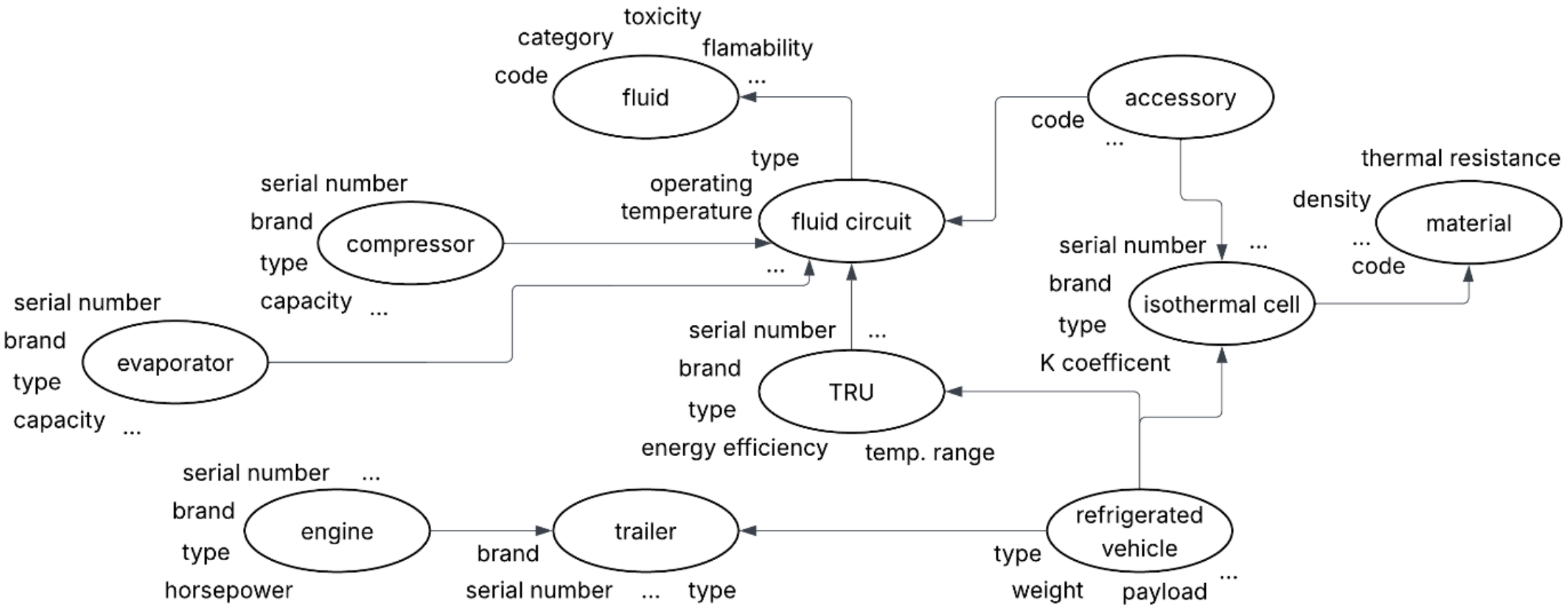

**Figure 2: Entity-Relationship diagram associated to a data model describing refrigerated vehicles**

### 3.3. Data quality assessment and correction (DQAC)

DQAC is crucial to analytics to ensure compliance with regulatory requirements, and because analytics built on bad quality data is known to carry a high financial cost and potentially lead to disastrous consequences (Redman, 2016; Haug et al., 2011). Furthermore, leveraging bad data can lead to ineffective resource allocation and misguided decisions, and thus undermine environmental impact reduction strategies.

Thinking about data in terms of entities facilitates the enumeration of many data quality and integrity constraints. Following our example, let us consider $F_C$ the set of constraints applicable to the entities of $C$, with $g \in F_C$: $interior_surface > 0$, and an entity $c \in C$, from Table 3, if $c[cell_id] = C4$ then $c[surface] = -21$ which is in direct violation of constraint $g$.

This violation would not have easily been identified in Tables 1 and 2, nor the constraint itself on Table 1. Data quality is mainly assessed in relation to real-world objects and processes and hence relies on the ability to uniquely identify a representation of real-life entities and their attributes. For example, business rules validation, one of the most used data quality assessment techniques, expressed data quality rules that data should follow such as "a refrigerated vehicle is necessarily associated to an insulated body".

Another critical example in the cold chain is duplicate detection, also referred to as entity resolution and record matching. In this context, duplicates are records that represent the same real-world element, and the objective is to produce a dataset of unique entities by identifying and merging such records. Comparisons are typically performed on key attributes, ideally corresponding to the entity's natural key, but may be extended to additional discriminative attributes when needed. For example, records C1 and C3 in Table 3 share identical values on key attributes $\{serial_number, brand\}$ and therefore represent the same insulated body
Duplicate detection is also essential for resolving other data quality problems, such as inconsistencies and missing values. In Table 3, records C1 and C3 differ in their reported interior surface values, although a real insulated body can have only one interior surface measure; resolving the duplication enables correction of this inconsistency. Similarly, missing values can be inferred from its duplicates. For example, if $c_1, c_2 \in C$ such that $c_1[cell_id] = C1$ and $c_2[cell_id] = C2$, and $c_1, c_2$ are duplicates, then $c_3[operational] \leftarrow c_1[operational]$.

This process is particularly important in cold chain data ecosystems, which often require the integration of multiple heterogeneous data sources. However, as it typically involves pairwise comparisons and exhibits quadratic complexity, it is computationally expensive and thus increases data processes' environmental footprint. Although techniques like blocking (Christen, 2012) can reduce the number of comparisons, the process remains resource intensive. In the simple previous example, resolving duplicates across a dataset of three records requires six comparisons, whereas a dataset of 100,000 records would entail nearly 15 billion comparisons. And although techniques such as blocking and filtering reduce the total number of comparisons, it remains a very computationally expensive problem that contributes to the environmental footprint of data ecosystems.

## 4. Cold-Chain Obstacles

Because most data technologies are designed to be domain-agnostic, they are not always well suited to the particularities of the sector. This section focuses on two major domain-specific challenges that complicate data engineering in the cold chain: the interdependence of entities and the impact on data modelling of multiple, sometimes conflicting, regulations and standards. While other constraints, such as the time sensitivity of data-driven decisions, also play an important role, they are beyond the scope of this article.

### 4.1. Entities interdependence

Operations in the cold chain rely in general on complex technological systems composed of layers of interdependent elements that can individually and collectively impact the system's overall performance and environmental footprint. Thus, to study a technology's performance, especially in the context of environmental impact mitigation, these components must also be modelled.

For example, many components from a refrigerated vehicle play a significant role in its energy consumption and thus environmental impact, as non-exhaustively represented in Figure 3. The insulated body and the transport refrigeration unit (TRU). TRU's efficiency is impacted by its compressor(s), and the body's insulation performance is impacted by the materials it is made of.

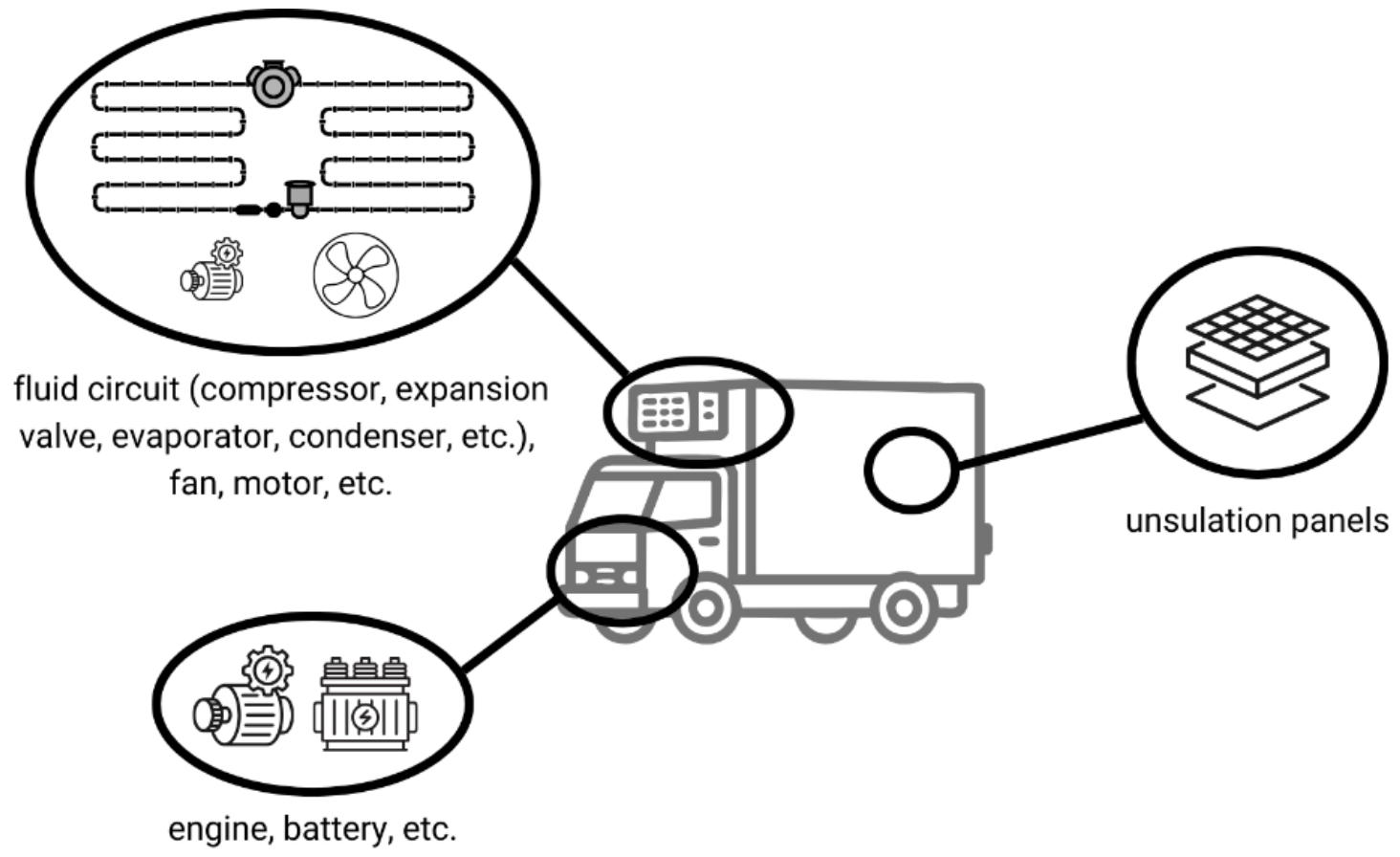

**Figure 3: Illustration of some technologies composing a refrigerated vehicle**

To assess the impact of such systems, by leveraging data from Figure 4, each unique combination of components must be treated as a distinct entity. For example, when focusing on the role played by insulated bodies, TRUs, and compressors (assuming for simplicity each vehicle contains a single TRU and each TRU a single compressor) in the environmental impact of a refrigerated vehicle, the uniqueness of a refrigerated vehicle depends on the combined uniqueness of its cell, TRU, and compressor.

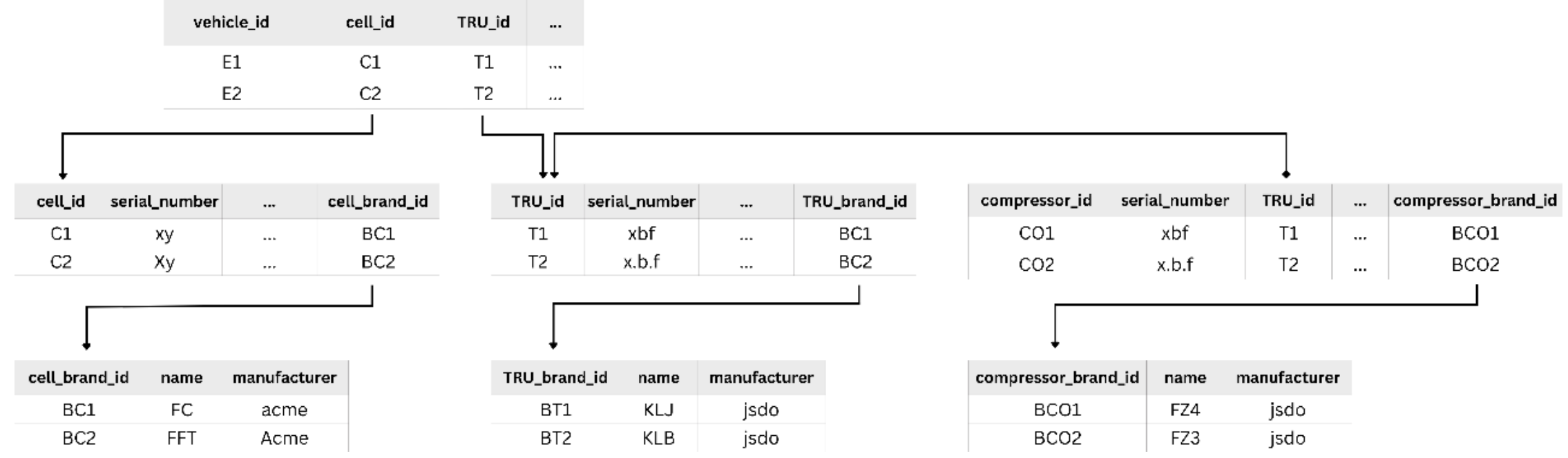

| vehicle_id | cell_id | TRU_id | ... |
|---|---|---|---|
| E1 | C1 | T1 | ... |
| E2 | C2 | T2 | ... |

| cell_id | serial_number | ... | cell_brand_id |
|---|---|---|---|
| C1 | xy | ... | BC1 |
| C2 | Xy | ... | BC2 |

| TRU_id | serial_number | ... | TRU_brand_id |
|---|---|---|---|
| T1 | xbf | ... | BC1 |
| T2 | x.b.f | ... | BC2 |

| compressor_id | serial_number | TRU_id | ... | compressor_brand_id |
|---|---|---|---|---|
| CO1 | xbf | T1 | ... | BCO1 |
| CO2 | x.b.f | T2 | ... | BCO2 |

| cell_brand_id | name | manufacturer |
|---|---|---|
| BC1 | FC | acme |
| BC2 | FFT | Acme |

| TRU_brand_id | name | manufacturer |
|---|---|---|
| BT1 | KLJ | jsdo |
| BT2 | KLB | jsdo |

| compressor_brand_id | name | manufacturer |
|---|---|---|
| BCO1 | FZ4 | jsdo |
| BCO2 | FZ3 | jsdo |

**Figure 4: data describing two refrigerated vehicles and their components**

Following collective entity resolution guidelines (Bhattacharya and Getoor, 2007), comparing vehicles records requires to compare the pair of each associated component records, as illustrated in Figure 5.
Considering that each component is uniquely identifiable by the combination of its brand and serial number, resolving this simple problem requires up to 7 pairwise comparisons. Most mainstream data techniques are not adapted to operate on such intricate interdependent entity structures, and performing duplicate detection on them is very computationally expensive. Moreover, data quality issues affecting any component's representation can propagate through the resolution process and compromise its outcome. While this phenomenon is only illustrated on entity resolution in this paper, similar limitations arise for other data engineering and science tasks, including rule validation, functional dependency discovery, master data management, and governance-related activities such as data lineage, security, and auditability.

Further complicating the application of data modelling in the cold chain, test results data and measurement data, are typically accompanied by uncertainties, which are critical for decision-making in safety and quality-sensitive cold chain operations. However, most classical data modelling and processing techniques are not designed to explicitly include such uncertainty-related metadata.

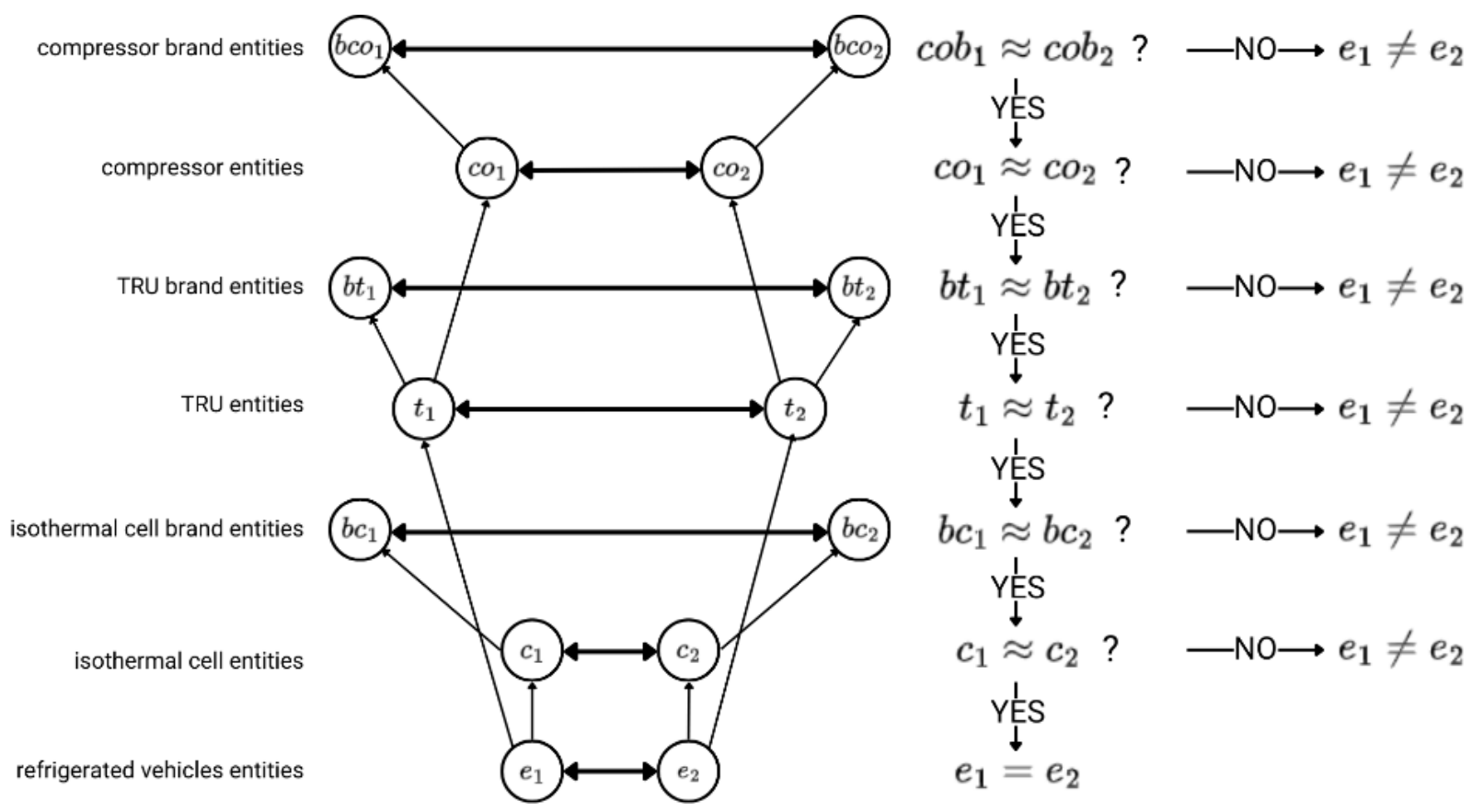

**Figure 5: Illustration of the collective entity resolution problem associated with refrigerated vehicles records**

### 4.2. Heterogeneous semantic references

The multitude of regulations and standards governing the cold chain introduces additional complexity into data modelling, as the definition of entities may vary across policies. Each regulatory framework may rely on its own semantic reference, sometimes partially conflicting with others. As a result, many data sources in the cold chain are modelled from the perspective of a specific regulation or operational context, complicating the integration of heterogeneous datasets, and the interoperability across systems.

For example, an application designed to ensure compliance with European lifting equipment inspection regulations may identify refrigerated trucks primarily through their trailers, with limited information on the insulated body. By contrast, one developed to monitor compliance with the ATP defines a refrigerated vehicle through its insulated body and therefore captures this information in detail. Integrating such data sources is challenging, as key attributes required for entity matching may be missing.

## 5. Foundations of the cold chain data management

To be effective, data processes should be adapted to their applications' contexts, although most data technologies were designed to be domain-agnostic (Goknil et al., 2023) and are not consistently aligned with domain-specificities. Certain sectors have created specialized branches of data management and science to more effectively meet their requirements, such as medical data management (Leiner et al., 2003).

To the authors' knowledge, there is no sub-field dedicated to data management in the cold chain; yet several domain specificities hinder the comprehensive adoption of data science in the sector. The authors advocate for the creation of this discipline to foster the development of adapted tools, provide guidance to the data professionals in the sector, and campaign for increased semantic harmonization in the sector. Such discipline would strengthen the sector's sustainability and enhance its performance, particularly as increased data and technological integration have been identified as one of the main drivers of evolution for the sector (Goknil et al., 2023; Alherimi and Ben-Daya, 2025).

This section discussed key concepts and orientations that, in the authors' view, are essential for the establishment of the *cold chain data management* discipline.

### 5.1. The data value chain

Although no single "one-size-fits-all" strategy can be universally applied across the cold chain sector, a common principle emerges: data science should be treated as a core activity rather than a supporting tool. From this perspective, data is considered a product, and data operations are structured as a value chain operating within a data ecosystem. This approach enables more efficient data management, aligning with the sector's goal of reducing its environmental footprint.

As described in section 2 of this paper, most existing studies have primarily focused on individual components of the data value chain. However, Marcy et al. (2024) started implementing this strategy through the development of REFTRAN (*REFrigerated TRansportation ANalytics*), a complete data ecosystem dedicated to refrigerated transportation.

### 5.2. Re-investing in data modelling and semantic references

Although data modelling is a long-standing discipline and a cornerstone of data management, it has increasingly been sidelined by developers, who often perceive it as overly time-consuming and instead rely on automated solutions such as object-relational mapping (ORM) frameworks. While these tools offer practical advantages, their widespread use can lead to a loss of expert knowledge and semantic precision and thus may result in misrepresentations of key real-world entities.

Accurately modelling entities in their full complexity remains a prerequisite for analyzing sector-specific technologies and operations and for quantifying and reducing their environmental impact. Reinvesting time and resources in data modelling and in the development of strong semantic references, such as ontologies, is therefore essential. Although this challenge extends beyond the cold chain sector, it is particularly acute in this domain due to the complexity of its operations and diversity of its data. Simultaneously, due to the sector's maturity and extensive domain expertise, it provides fertile conditions for the development and implementation of rigorous semantic modeling approaches.

Furthermore, to avoid conflicting semantics, data should be modelled from a perspective that is common to all present and future relevant contexts, such as the intrinsic definitions of real-world elements, rather than the constraints of a specific activity. Indeed, modelling a refrigerated vehicle as a modular composition of entities (insulated body, trailer, cooling systems, etc.), is compatible with multiple regulatory and operational viewpoints without compromising data coherence.

### 5.3. Develop tools or adapt existing ones

Data ecosystems should be designed with full consideration of domain-specific constraints and specificities, which may require the use of sector-specific data tools or the adaptation of mainstream solutions. For example, Marcy et al. (2025) propose RED2Hunt, a human-in-the-loop framework for improving duplicate detection and resolving data inconsistencies when identification of unique entities is interdependent between entity classes. It leverages domain experts' knowledge, and artificial elements within relational databases (which are normally excluded from DQAC techniques), to considerably improve the problem's resolution and reduce its complexity, and hence its environmental impact.

A limited shift toward adapting data cleaning and data sciences applications to Industry 4.0 contexts has been observed, with a focus on IoT and CPS (Ding et al., 2022; Goknil et al., 2023), but without specialization for the cold chain.

Building on this, the authors advocate for the development of human-in-the-loop data tools in *cold chain data management*. While a prevailing trend in data science aims to minimize reliance on domain experts, often perceived as a scarce resource, the cold chain sector benefits from a wealth of expertise. As such, domain knowledge should be actively leveraged as a central component of effective data ecosystem design.

## 6. Conclusion

While the adoption of new data technologies offers amazing opportunities to the cold chain sector to control and reduce their environmental footprint, they could also become an environmental burden themselves. To balance this tradeoff in favor of more sustainability, it is important for companies within the sector to consider data as a product and data science as a value chain to optimize, and to design their data ecosystem with performance requirements in mind.

Several domain-specific constraints complexify the deployment of data sciences within the cold chain industry, such as the high level of granularity that characterizes modelling of sector data. The authors argue in favor of the creation of a *cold chain data management* discipline. This discipline would gather researchers and practitioners from both data and the cold chain sector to develop data tools optimized to the reality of the domain, lobby for semantic harmonization within the sector, train cold chain data professionals.

The authors also believe that such discipline could be a source of positive evolution for data science as an influence for an improved uncertainty assessment and more regulations. Especially without an efficient and robust cold chain, data centers could not operate consistently, and thus the democratization of data science would be at risk.

## REFERENCES

Alherimi, N., Ben-Daya, M., 2025. A Systematic Review on the Intersection of the Cold Chain and Digital Transformation. Sustainability, 17(24), 11202.

Arora, M., Ahmad, V., Kumar, R., Yamsani, N., Amir, M., Arora, J., 2024. Blockchain and Internet-of-Things: A Technological Solution for Issues in Cold Chain. Proceeding of the 3rd International Conference on Sentiment Analysis and Deep Learning. ICSADL IEEE, 632-636.

Bhagat, P. R., Naz, F., Magda, R., 2022. Role of Industry 4.0 Technologies in enhancing sustainable firm performance and green practices. Acta Polytechnica Hungarica, 19(8), 229-248.

Bhattacharya, I., Getoor, L., 2007. Collective entity resolution in relational data. ACM Transactions on Knowledge Discovery from Data, TKDD, 1(1), 5-es.

Bottani, E., Casella, G., Nobili, M., Tebaldi, L., 2022. An analytic model for estimating the economic and environmental impact of food cold supply chain. Sustainability, 14(8), 4771.

Capo, C., Le Guilly M, Pape L, Scuturici M, Petit, J.M., Revellin, R., Bonjour, J., Cavalier, G., 2020. Ageing of refrigerated transport vehicles: development of a numerical predictive model. Proceeding of the 6th IIR Conference on Sustainability and the Cold Chain, Nantes, France.

Chaudhuri, A., Dukovska-Popovska, I., Subramanian, N., Chan, H. K., Bai, R., 2018. Decision-making in cold chain logistics using data analytics: a literature review. The International Journal of Logistics Management, 29(3), 839-861.

Christen P., 2012. Data Matching: Concepts and Techniques for Record Linkage, Entity Resolution, and Duplicate Detection, Springer, Berlin, Germany. 272p.

Corbett, C.J., 2018. How Sustainable Is Big Data? Journal Production and Operations Management, 27(9) ISSN 1059-1478.Ding, X., Wang, H., Li, G., Li, H., Li, Y., Liu, Y., 2022. IoT data cleaning techniques: A survey. Intelligent and Converged Networks, 3(4), 325-339.

Draisbach U., Christen P., Naumann F., 2020. Transforming pairwise duplicates to entity clusters for high-quality duplicate detection, ACM J. Data Inf. Qual., vol. 12, no. 1, Art. 3.

Fabris F., Artuso P., Marinetti S., Minetto S., Rossetti A., 2022. Cooling unit impact on energy and emissions of a refrigerated light truck. Applied Thermal Engineering, 216, 119132.

Fatorachian, H., Pawar, K., 2025. Waste efficiency in cold supply chains through industry 4.0-enabled digitalisation. International Journal of Sustainable Engineering, 18(1), 2461564.

Gillespie, J., da Costa, T. P., Cama-Moncunill, X., Cadden, T., Condell, J., Cowderoy, T., Ramanathan, R., 2023. Real-time anomaly detection in cold chain transportation using IoT technology. Sustainability, 15(3), 2255.

Goknil, A., Nguyen, P., Sen, S., Politaki, D., Niavis, H., Pedersen, K. J., Ziegenbein, A., 2023. A systematic review of data quality in CPS and IoT for industry 4.0. ACM Computing Surveys, 55(14s), 1-38.

Haug, A., Zachariassen, F., Van Liempd, D., 2011. The costs of poor data quality. Journal of Industrial Engineering and Management, JIEM, 4(2), 168-193.

Kale, S. D., 2022. Predictive analytics for cold chain break detection. Preprint available at Research Square.

Kavididevi, V., Monikapreethi, S. K., Rajapriya, M., Juliet, P. S., Yuvaraj, S., Muthulekshmi, M., 2024. IoT-Enabled Reinforcement Learning for Enhanced Cold Chain Logistics Performance in Refrigerated Transport. Proceeding of the 2nd International Conference on Sustainable Computing and Smart Systems, IEEE ICSCSS, 379-384.

Khanna, A., Jain, S., Sah, A., Dangi, S., Sharma, A., Tiang, S. S., Lim, W. H., 2025. Generative AI and Blockchain-Integrated Multi-Agent Framework for Resilient and Sustainable Fruit Cold-Chain Logistics. Foods, 14(17), 3004.

Leiner, F., Gaus, W., Haux, R., Knaup-Gregori, P., 2003. Medical data management: a practical guide, Springer, USA: New York.

Lorenc, A., Czuba, M., Szarata, J., 2021. Big Data Analytics and Anomaly Prediction in the Cold Chain to Supply Chain Resilience. FME transactions, 49(2).

Li, C., Huang, M., 2023. Environmental Sustainability in the Age of Big Data: Opportunities and Challenges for Business and Industry. Environmental Science and Pollution Research **30**, 119001–119015.

Marcy M., Fertel C., Petit J.-M., Scuturici V.-M., Cavalier G., Bonjour J., 2024. Modelling and prediction of refrigerated trucks' environmental performance using real-life data and data science. Proceeding of the 8th IIR International Conference on Sustainability and the Cold Chain, Tokyo, Japan.

Marcy M., Petit J.-M., Scuturici M. , Bonjour J. , Fertel C., Cavalier G., 2025. Can Surrogate Keys Negatively Impact Data Quality? Proceedings of the 2025 Very Large Databases Conference, VLDB, London, UK.

Mohan, M., Amin, S., 2025. Green Cold Chain Logistics: Minimising Greenhouse Gas Emissions of Fresh Food Products in Transport Refrigeration Units. Logistics, 9(3), 112.

[Mtibaa, A., 2025. Refrigerant Leak Detection in Industrial Vapor Compression Refrigeration Systems. Doctoral dissertation, Université Paris sciences et lettres. 234p.

Redman, T. C., 2016. Bad data costs the US $3 trillion per year. Harvard Business Review, 22, 11-18.

Sartal, A., Bellas, R., Mejías, A. M., García-Collado, A., 2020. The sustainable manufacturing concept, evolution and opportunities within Industry 4.0: A literature review. Advances in Mechanical Engineering, 12(5).

Vrchota, J., Pech, M., Rolinek, L., Bednář, J., 2020. Sustainability outcomes of green processes in relation to industry 4.0 in manufacturing: Systematic review. Sustainability, 12(15), 5968.

Whitehead, B., Andrews, D., Shah, A., Maidment, G., 2014. Assessing the environmental impact of data centres part 1: Background, energy use and metrics. Building and Environment, 82, 151-159.

Wu, W., Shen, L., Zhao, Z., Harish, A. R., Zhong, R. Y., Huang, G. Q., 2023. Internet of everything and digital twin enabled service platform for cold chain logistics. Journal of Industrial Information Integration, 33, 100443.

Yao Y., Shi L., He J., Tian H., Wang X., Zhang X., Sun X., Lu B., Shu G., 2023. Dynamic analysis of refrigerated truck integrated with combined cooling and power cycle under various driving conditions. Energy Conversion and Management 278, 116716.

Zhao, B., Gui, H., Li, H., Xue, J., 2020. Cold chain logistics path optimization via improved multi-objective ant colony algorithm. Ieee Access, 8, 142977-142995.